\documentclass[10pt]{article}
\usepackage{mathrsfs}
\usepackage{amssymb}
\usepackage{amsmath}
\usepackage[all]{xy}

\def\Im{\mathop{\rm Im}\nolimits}
\def\Ker{\mathop{\rm Ker}\nolimits}

\def\mod{\mathop{\rm mod}\nolimits}
\def\Mod{\mathop{\rm Mod}\nolimits}
\def\Hom{\mathop{\rm Hom}\nolimits}

\def\Ext{\mathop{\rm Ext}\nolimits}

\def\add{\mathop{\rm add}\nolimits}

\def\Proj{\mathop{\rm Proj}\nolimits}
\def\proj{\mathop{\rm proj}\nolimits}
\def\Inj{\mathop{\rm Inj}\nolimits}
\def\inj{\mathop{\rm inj}\nolimits}

\def\inf{\mathop{\rm inf}\nolimits}
\def\sup{\mathop{\rm sup}\nolimits}
\def\dim{\mathop{\rm dim}\nolimits}

\def\Im{\mathop{\rm Im}\nolimits}
\def\Ker{\mathop{\rm Ker}\nolimits}

\def\mod{\mathop{\rm mod}\nolimits}
\def\Mod{\mathop{\rm Mod}\nolimits}

\def\id{\mathop{\rm id}\nolimits}

\def\sup{\mathop{\rm sup}\nolimits}
\def\inf{\mathop{\rm inf}\nolimits}
\def\add{\mathop{\rm add}\nolimits}

\def\G{\mathop{\mathscr{G}}\nolimits}

\def\dim{\mathop{\rm dim}\nolimits}

\def\Hom{\mathop{\rm Hom}\nolimits}
\def\Ext{\mathop{\rm Ext}\nolimits}
\def\sup{\mathop{\rm sup}\nolimits}
\def\lim{\mathop{\underrightarrow{\rm lim}}\nolimits}

\def\Con{\mathop{\rm Con}\nolimits}

\title{\large \bf Relative Singularity Categories
}
\author{Huanhuan Li\thanks{{\it E-mail address}: lihuanhuan0416@163.com}
\ and Zhaoyong Huang\thanks{{\it E-mail address}: huangzy@nju.edu.cn}\\
{\footnotesize \it Department of Mathematics, Nanjing University,
Nanjing 210093, Jiangsu Province, P.R. China}}
\date{}
\begin{document}
\baselineskip=16pt
\maketitle

\begin{abstract}We study the properties of the relative derived category $D_{\mathscr{C}}^{b}$($\mathscr{A}$) of an abelian
category $\mathscr{A}$ relative to a full and additive subcategory
$\mathscr{C}$. In particular, when $\mathscr{A}=A{\text -}\mod$ for a finite-dimensional algebra $A$ over
a field and $\mathscr{C}$ is a contravariantly finite subcategory of $A$-$\mod$ which is admissible and closed under direct summands,
the $\mathscr{C}$-singularity category $D_{\mathscr{C}{\text-sg}}$($\mathscr{A}$)=$D_{\mathscr{C}}^{b}$($\mathscr{A}$)/$K^{b}(\mathscr{C})$
is studied. We give a sufficient condition when this category is triangulated equivalent to the stable category of the Gorenstein category
$\mathscr{G}(\mathscr{C})$ of $\mathscr{C}$.
\end{abstract}
\vspace{0.5cm}

\centerline{\bf  1. Introduction}

\vspace{0.2cm} Let $A$ be a finite-dimensional algebra over a field. We denote by $A$-$\mod$ the category
of finitely generated left $A$-modules, and $A$-$\proj$ (resp. $A$-$\inj$) the full subcategory of $A$-$\mod$
consisting of projective (resp. injective) modules. We use
$K^b(A)$ and $D^b(A)$ to denote the bounded homotopy and derived categories of $A$-$\mod$ respectively,
and $K^b(A$-$\proj)$ (resp. $K^b(A$-$\inj)$) to denote the bounded homotopy category of $A$-$\proj$ (resp. $A$-$\inj$).

The composition functor $K^b(A$-$\proj) \to K^b(A)\to D^b(A)$ with the former functor the inclusion
functor and the latter one the quotient functor is naturally a fully faithful triangle functor, and then one
can view $K^b(A$-$\proj)$ as a triangulated subcategory of $D^b(A)$. In fact it is a thick one by [Bu, Lemma 1.2.1].
Consider the quotient
triangulated category $D_{sg}(A):=D^b(A)/K^b(A$-$\proj)$, which is the so-called ``singularity category". This category
was first introduced and studied  by Buchweitz in [Bu] where $A$ is assumed to be a left and right noetherian ring.
Later on Rickard proved in [R] that for a self-injective algebra $A$, this category is triangle-equivalent to the
stable category of $A$-$\mod$. This result was generalized to Gorenstein algebra by Happel in [H2]. Since
$A$ has finite global dimension  if and only if $D_{sg}(A)=0$, from this viewpoint $D_{sg}(A)$ measures the homological
singularity of the algebra $A$, we call it the singularity category after [O].

Besides, other quotient triangulated categories have been studied by many authors. Beligiannis considered the quotient
triangulated categories $D^b(R$-$\Mod)/K^b(R$-$\Proj)$ and $D^b(R$-$\Mod)/K^b(R$-$\Inj)$ for arbitrary ring $R$, where $R$-$\Mod$
is the category of left $R$-modules and $R$-$\Proj$ (resp. $R$-$\Inj$) is the full subcategory of $R$-$\Mod$
consisting of projective (resp. injective) modules (see [Be]). Let $\mathscr{A}$ be an abelian category.
A full and additive subcategory $\omega$ of $\mathscr{A}$ is called {\it self-orthogonal} if $\Ext_{\mathscr{A}}^i(M,N)=0$ for
any $M,N\in\omega$ and $i\geq 1$; in particular, an object $T$ in $\mathscr{A}$ is called {\it self-orthogonal} if $\Ext_{\mathscr{A}}^i(T,T)=0$
for any $i\geq 1$. Chen and Zhang studied in [CZ] the quotient triangulated category $D^b(A)/K^b(\add_A T)$ for a finite-dimensional
algebra $A$ and a self-orthogonal module $T$ in $A$-$\mod$, where $\add_A T$ is the full subcategory of
$A$-$\mod$ consisting of direct summands of finite direct sums of $T$. Recently Chen studied in [C2] the relative singularity
category $D_\omega(\mathscr{A}):=D^b(\mathscr{A})/K^b(\omega)$ for an arbitrary abelian category $\mathscr{A}$ and an arbitrary
self-orthogonal, full and additive subcategory $\omega$ of $\mathscr{A}$.

For an abelian category $\mathscr{A}$ with enough projective objects, the Gorenstein derived category
$D^*_{gp}(\mathscr{A})$ of $\mathscr{A}$ was introduced by Gao and Zhang in [GZ], where $* \in \{{\rm blank},-,b\}$.
It can be viewed as a generalization of the usual derived category $D^*(\mathscr{A})$ by using Gorenstein projective
objects instead of projective objects and $\mathscr{G}\mathscr{P}$-quasi-isomorphisms instead of quasi-isomorphisms,
where $\mathscr{G}\mathscr{P}$ means ``Gorenstein projective". For Gorenstein projective modules and Gorenstein projective
objects, we refer to [AuB], [EJ1], [EJ2], [Ho] and [SSW]. Asadollahi, Hafezi and Vahed studied
in [AHV] the relative derived category $D^*_{\mathscr{C}}(\mathscr{A})$ for an arbitrary abelian category $\mathscr{A}$
with respect to a contravariantly finite subcategory $\mathscr{C}$, where $* \in \{{\rm blank},-,b\}$, and they
pointed out that $K^{b}(\mathscr{C}$) can be viewed as a triangulated subcategory of
$D^b_{\mathscr{C}}(\mathscr{A})$.

Given a finite-dimensional algebra $A$ over a field and a full and additive subcategory $\mathscr{C}$ of $\mathscr{A}(=A$-$\mod)$ closed under
direct summands, it follows from [BD] that $K^{b}(\mathscr{C}$) is a Krull-Schmidt category and hence can be viewed as a
thick triangulated subcategory of $D^b_{\mathscr{C}}(\mathscr{A})$. If the quotient triangulated category
$D_{\mathscr{C}{\text-}sg}(\mathscr{A}):=D^b_{\mathscr{C}}(\mathscr{A})/K^{b}(\mathscr{C})$ is considered, then it is natural to ask
whether $D_{\mathscr{C}{\text-}sg}(\mathscr{A})$ shares some nice properties of $D_{sg}(A)$. The aim of this paper is to study this question.

In Section 2, we give some terminology and some preliminary results.

In Section 3, for an abelian category $\mathscr{A}$ and a full and additive subcategory $\mathscr{C}$ of $\mathscr{A}$,
we prove that if $\mathscr{C}$ is admissible, then the composition functor
$\mathscr{A}\to K^{b}(\mathscr{A})\to D_{\mathscr{C}}^{b}(\mathscr{A})$ is fully faithful,
where the former functor is the inclusion functor and the latter one is the quotient functor.
Let $\mathscr{C}$ be a contravariantly finite subcategory of $\mathscr{A}$ and $\mathscr{D}\subseteq\mathscr{A}$ a
subclass of $\mathscr{A}$. We introduce a dimension denoted by $\mathscr{C}\mathscr{D}$-$\dim M$ which is called the
{\it $\mathscr{C}$-proper $\mathscr{D}$-dimension} of an object $M$ in $\mathscr{A}$. By choosing a left $\mathscr{C}$-resolution
$C^\bullet_M$ of $M$, we get a functor $\Ext_{\mathscr{C}}^n(M,-):=H^{n}\Hom_\mathscr{A}(C^{\bullet}_M,-)$ for any $n \in \mathbb{Z}$.
Then by using the properties of this functor we obtain some equivalent characterizations for $\mathscr{C}\mathscr{C}$-$\dim M$ being finite.

In Section 4, we introduce the $\mathscr{C}$-singularity category $D_{\mathscr{C}{\text-}sg}(\mathscr{A}):=
D^{b}_\mathscr{C}(\mathscr{A})$ /$K^{b}(\mathscr{C})$, where $\mathscr{A}=A$-$\mod$ and $\mathscr{C}$ is a contravariantly
finite, full and additive subcategory of $\mathscr{A}$ which is admissible and closed under direct summands. We prove that if
$\mathscr{C}\mathscr{C}$-$\dim \mathscr{A}<\infty$, then $D_{\mathscr{C}{\text-}sg}(\mathscr{A})=0$.
As a consequence, we get that if $A$ is of finite representation type, then $\mathscr{C}\mathscr{C}$-dim$\mathscr{A}<\infty$
if and only if $D_{\mathscr{C}{\text-}sg}(\mathscr{A})=0$. Let $\mathscr{G}(\mathscr{C})$ be the Gorenstein category
of $\mathscr{C}$ and $\varepsilon$ the collection of all $\Hom_\mathscr{A}(\mathscr{C},-)$-exact complexes of the form:
$0\to L\to M\to N\to 0$ with $L, M, N\in \mathscr{G}(\mathscr{C})$. By [B\"{u}] (or [Q])
$(\mathscr{G}(\mathscr{C}),\varepsilon)$ is an exact category; moreover, it is a Frobenius category with $\mathscr{C}$
the subcategory of projective-injective objects, see [H1]. We prove that if $\mathscr{C}\mathscr{G}(\mathscr{C})$-dim $\mathscr{A}<\infty$,
then the natural functor $\theta: \mathscr{G}(\mathscr{C})\to D_{\mathscr{C}{\text-}sg}(\mathscr{A})$ induces a triangle-equivalence
$\theta^{'}: \underline{\mathscr{G}(\mathscr{C})} \to D_{\mathscr{C}{\text-}sg}(\mathscr{A})$, where $\underline{\mathscr{G}(\mathscr{C})}$
is the stable category of $\mathscr{G}(\mathscr{C})$.

\vspace{0.5cm}

\centerline{\bf 2. Preliminaries}

\vspace{0.2cm}
Throughout this paper, $\mathscr{A}$ is an abelian category, $C(\mathscr{A})$ is the category of complexes of objects
in $\mathscr{A}$, $K^{*}(\mathscr{A})$ is the homotopy category of $\mathscr{A}$ and $D^{*}(\mathscr{A})$ is the
usual derived category by inverting the quasi-isomorphisms in $K^{*}(\mathscr{A})$, where $* \in \{{\rm blank},-,b\}$.
We will use the formula $\Hom_{K(\mathscr{A})}(X^{\bullet},Y^{\bullet}[n])$=$H^n\Hom_\mathscr{A}(X^{\bullet},Y^{\bullet})$
for any $X^{\bullet},Y^{\bullet}\in C(\mathscr{A})$ and $n\in \mathbb{Z}$ (the ring of integers).

Let $$X^{\bullet}:=\cdots \longrightarrow X^{n-1}\buildrel {d^{n-1}_{X}} \over \longrightarrow X^{n}\buildrel {d^{n}_{X}}
\over\longrightarrow X^{n+1} \to \cdots$$ be a complex and $f: X^{\bullet}\to Y^{\bullet}$ a cochain map in $C(\mathscr{A})$.
Recall that $X^{\bullet}$ is called {\it acyclic} (or {\it exact}) if $H^{i}(X^{\bullet})=0$ for any $i\in \mathbb{Z}$, and $f$
is called a {\it quasi-isomorphism} if $H^{i}(f)$ is an isomorphism for any $i\in \mathbb{Z}$.

From now on, we fix a full and additive subcategory $\mathscr{C}$ of $\mathscr{A}$.

\vspace{0.2cm}
{\bf Definition 2.1.}  Let $X^{\bullet},Y^{\bullet}$ and $f$ be as above.

(1) ([EJ2]) $X^{\bullet}$ in $C(\mathscr{A})$ is called {\it $\mathscr{C}$-acyclic} or {\it $\Hom_\mathscr{A}(\mathscr{C},-)$-exact}
if the complex $\Hom_{\mathscr{A}}(C,X^{\bullet})$ is acyclic for any $C\in\mathscr{C}$. Dually, a {\it $\Hom_{\mathscr{A}}(-,\mathscr{C})$-exact
complex} is defined.

(2) $f$ is called a {\it $\mathscr{C}$-quasi-isomorphism} if the cochain map $\Hom_\mathscr{A}(C,f)$ is a quasi-isomorphism for any $C\in \mathscr{C}$.

\vspace{0.2cm}
{\bf Remark 2.2.} (1) We use $\Con(f)$ to denote the mapping cone of $f: X^{\bullet}\to Y^{\bullet}$. It is well known that $f$
is a quasi-isomorphism if and only if $\Con(f)$ is acyclic; analogously, $f$ is a $\mathscr{C}$-quasi-isomorphism if and only if
$\Con(f)$ is $\mathscr{C}$-acyclic.

(2) We use $\mathscr{P(A)}$ to denote the full subcategory of $\mathscr{A}$ consisting of projective objects.
If $\mathscr{A}$ has enough projective objects,
then every quasi-isomorphism is a $\mathscr{P(A)}$-quasi-isomorphism; and if $\mathscr{P(A)}\subseteq \mathscr{C}$,
then every $\mathscr{C}$-quasi-isomorphism is a quasi-isomorphism.

\vspace{0.2cm}
We use $K_{ac}^{*}(\mathscr{A})$ (resp. $K_{\mathscr{C}{\text -ac}}^{*}(\mathscr{A})$) to denote
the full subcategory of $K^{*}(\mathscr{A})$ consists of acyclic complexes (resp. $\mathscr{C}$-acyclic complexes).

\vspace{0.2cm}

{\bf Lemma 2.3.} {\it Let $X^{\bullet}$ be a complex in $C(\mathscr{A})$. Then $X^{\bullet}$ is $\mathscr{C}$-acyclic
if and only if the complex $\Hom_\mathscr{A}(C^{\bullet},X^{\bullet})$ is acyclic for any $C^{\bullet}\in K^{-}(\mathscr{C})$.}

\vspace{0.2cm}

{\it Proof.} See [CFH, Lemma 2.4]. \hfill{$\square$}

\vspace{0.2cm}

{\bf Lemma 2.4.} {\it (1) Let $C^{\bullet}$ be a complex in $K^{-}(\mathscr{C})$ and f : $X^{\bullet}\to C^{\bullet}$
a $\mathscr{C}$-quasi-isomorphism in $C(\mathscr{A})$. Then there exists a cochain map $g: C^{\bullet}\to X^{\bullet}$
such that $fg$ is homotopic to $\id_{C^{\bullet}}$.

(2) Any $\mathscr{C}$-quasi-isomorphism between two complexes in $K^{-}(\mathscr{C})$ is a homotopy equivalence.}

\vspace{0.2cm}
{\it Proof.} (1) Consider the distinguished triangle:
$$X^{\bullet}\buildrel {f} \over \longrightarrow C^{\bullet}\to
\Con(f)\to X^{\bullet}[1]$$ in $K(\mathscr{A})$ with $\Con(f)$ $\mathscr{C}$-acyclic. By applying the functor
$\Hom_{K(\mathscr{A})}(C^\bullet,-)$ to it, we get an exact sequence:
$$\Hom_{K(\mathscr{A})}(C^\bullet,X^\bullet)\buildrel {\Hom_{K(\mathscr{A})}(C^\bullet,f)} \over \longrightarrow
\Hom_{K(\mathscr{A})}(C^\bullet,C^\bullet)\to \Hom_{K(\mathscr{A})}(C^\bullet,\Con(f)).$$
It follows from Lemma 2.3 that $\Hom_{K(\mathscr{A})}(C^\bullet,\Con(f))\cong H^0\Hom_\mathscr{A}(C^{\bullet},\Con(f))=0$.
So there exists a cochain map $g: C^{\bullet}\to X^{\bullet}$ such that $fg$ is homotopic to $\id_{C^{\bullet}}$.

(2) Let $f :X^{\bullet}\to Y^{\bullet}$ be
a $\mathscr{C}$-quasi-isomorphism with $X^{\bullet}, Y^{\bullet}$ in $K^{-}(\mathscr{C})$. By (1), there exists a
cochain map $g: Y^{\bullet}\to X^{\bullet}$, such that $fg$ is homotopic to $\id_{Y^{\bullet}}$. By (1) again,
there exists a cochain map $g': X^{\bullet}\to Y^{\bullet}$, such that $gg'$ is homotopic to $\id_{X^{\bullet}}$.
Thus $f=g'$ in $K(\mathscr{A})$ is a homotopy equivalence. \hfill{$\square$}

\vspace{0.2cm}

{\bf Definition 2.5.} (1) ([AuR]) Let $\mathscr{C}\subseteq\mathscr{D}$ be
subcategories of $\mathscr{A}$. The morphism $f: C\to D$ in
$\mathscr{A}$ with $C\in\mathscr{C}$ and $D\in\mathscr{D}$ is called a {\it right $\mathscr{C}$-approximation} of $D$
if for any morphism $g: C^{'} \to D$ in $\mathscr{A}$ with $C^{'}\in\mathscr{C}$, there exists a morphism $h: C^{'}\to C$
such that the following diagram commutes:
$$\xymatrix{ & C^{'} \ar[d]^{g} \ar@{-->}[ld]_{h}\\
C \ar[r]^{f} & D.}$$
If each object in $\mathscr{D}$ has a right $\mathscr{C}$-approximation,
then $\mathscr{C}$ is called {\it contravariantly finite} in $\mathscr{D}$.

(2) ([C1]) A contravariantly finite subcategory $\mathscr{C}$ of $\mathscr{A}$ is called {\it admissible} if any right $\mathscr{C}$-approximation is epic.
In this case, every $\mathscr{C}$-acyclic complex is acyclic.

\vspace{0.2cm}

The following definition is cited from [B\"{u}], see also [Q] and [K].

\vspace{0.2cm}

{\bf Definition 2.6.} Let $\mathscr{B}$ be an additive category. A {\it kernel-cokernel pair} $(i,p)$ in $\mathscr{B}$ is a pair
of composable morphisms $L\buildrel {i}\over \longrightarrow M\buildrel {p} \over \longrightarrow N$ such that
$i$ is a kernel of $p$ and $p$ is a cokernel of $i$. If a class $\varepsilon$ of kernel-cokernel pairs on $\mathscr{B}$
is fixed, an {\it  admissible monic} (sometimes called {\it inflation}) is a morphism $i$ for which there exists a morphism $p$
such that $(i,p)\in \varepsilon$. {\it Admissible epics} (sometimes called {\it deflations}) are defined dually.

{\it An exact category} is a pair $(\mathscr{B},\varepsilon)$ consisting of an additive category $\mathscr{B}$
and a class of kernel-cokernel pairs $\varepsilon$ on $\mathscr{B}$ with $\varepsilon$ closed under isomorphisms satisfying the following axioms:

[E0] For any object $B$ in $\mathscr{B}$, the identity morphism $\id_B$ is both an admissible monic and an admissible epic.

[E1] The class of admissible monics is closed under compositions.

[${\text E1^{op}}$] The class of admissible epics is closed under compositions.

[E2] The push-out of an admissible monic along an arbitrary morphism exists and yields an admissible monic.

[${\text E2^{op}}$] The pull-back of an admissible epic along an arbitrary morphism exists and yields an admissible epic.

Elements of $\varepsilon$ are called {\it short exact sequences} (or {\it conflations}).

\vspace{0.2cm}

Let $\mathscr{B}$ be a triangulated subcategory of a triangulated category $\mathscr{K}$ and
$S$ the compatible multiplicative system determined by $\mathscr{B}$. In the Verdier quotient category $\mathscr{K}/\mathscr{B}$,
each morphism $f: X\to Y$ is given by an equivalence class of right
fractions $f/s$ or left fractions $s\backslash f$ as presented by $X\buildrel {s}\over \Longleftarrow Z\buildrel {f}\over \longrightarrow Y$
or $X\buildrel {f}\over \longrightarrow Z\buildrel {s}\over \Longleftarrow Y$, where the doubled arrow means $s\in S$.

\vspace{0.5cm}

\centerline{\bf 3. $\mathscr{C}$-derived categories}

\vspace{0.2cm}

For a subclass $\mathscr{C}$ of objects in a triangulated category $\mathscr{K}$, it is known that the full subcategory
$\mathscr{C}^\perp=\{X\in\mathscr{K}\mid \Hom_\mathscr{K}(C[n],X)=0$ for any $C\in \mathscr{C}$ and $n\in \mathbb{Z}\}$
is a triangulated subcategory of $\mathscr{K}$ and is closed under direct summands, and hence is thick ([R]).
It follows that $K_{\mathscr{C}{\text -ac}}^{*}(\mathscr{A})$ is a thick subcategory of $K^{*}(\mathscr{A})$.

\vspace{0.2cm}

{\bf Definition 3.1.} ([V]) The Verdier quotient category $D_{\mathscr{C}}^{*}$($\mathscr{A}):
=K^{*}(\mathscr{A})/K_{\mathscr{C}{\text -ac}}^{*}(\mathscr{A})$
is called the {\it $\mathscr{C}$-derived category} of $\mathscr{A}$, where $*\in\{{\rm blank},-,b\}$.

\vspace{0.2cm}

{\bf Example 3.2.} (1) If $\mathscr{A}$ has enough projective objects and $\mathscr{C}=\mathscr{P(A)}$,
then $D_{\mathscr{C}}^{*}$($\mathscr{A}$) is the usual derived category $D^{*}(\mathscr{A})$.

(2) If $\mathscr{A}$ has enough projective objects and $\mathscr{C}=\mathscr{G(A)}$ (the full subcategory of $\mathscr{A}$ consisting of
Gorenstein projective objects), then $D_{\mathscr{C}}^{*}$($\mathscr{A}$) is the Gorenstein
derived category $D_{gp}^{*}$($\mathscr{A}$) defined in [GZ].

(3) Let $R$ be a ring and $\mathscr{A}=R{\text-} \Mod$. If $\mathscr{C}=\mathscr{P}\mathscr{P}(R)$
(the full subcategory of $R{\text-} \Mod$ consisting of pure projective modules), then $D_{\mathscr{C}}^{*}$($\mathscr{A}$)
is the pure derived category $D_{pur}^{*}$($\mathscr{A}$) in [ZH].

\vspace{0.2cm}

{\bf Proposition 3.3.} {\it ([AHV]) (1) $D_{\mathscr{C}}^{-}(\mathscr{A})$ is a triangulated subcategory of
$D_{\mathscr{C}}(\mathscr{A}$), and $D_{\mathscr{C}}^{b}(\mathscr{A})$ is a triangulated subcategory
of $D_{\mathscr{C}}^{-}(\mathscr{A})$.

(2) For any $C^{\bullet}\in K^{-}(\mathscr{C})$ and $X^{\bullet}\in C(\mathscr{A})$, there exists an isomorphism of abelian groups:
$$\Hom_{K(\mathscr{A})}(C^{\bullet},X^{\bullet}) \cong  \Hom_{D_\mathscr{C}(\mathscr{A})}(C^{\bullet},X^{\bullet}).$$

(3) Let $\mathscr{C}\subseteq\mathscr{A}$ be admissible. Then the composition functor
$\mathscr{A}\to K^{b}(\mathscr{A})\to D_{\mathscr{C}}^{b}(\mathscr{A})$ is fully faithful,
where the former functor is the inclusion functor and the latter one is the quotient functor.}

\vspace{0.2cm}
{\it Proof.} In the following,
each morphism in $D_{\mathscr{C}}^{*}(\mathscr{A})$ will be denoted by the equivalence class of right fractions,
where $*\in\{{\rm blank},-,b\}$.

(1) We only prove the first assertion, the second one can be proved similarly.

Note that $D_{\mathscr{C}}^{-}(\mathscr{A})=K^{-}(\mathscr{A})/K^{-}(\mathscr{A})\bigcap K_{\mathscr{C}{\text -ac}}(\mathscr{A})$
and $D_{\mathscr{C}}$($\mathscr{A})=K(\mathscr{A})/K_{\mathscr{C}{\text -ac}}(\mathscr{A})$. By [GM, Proposition 3.2.10],
it suffices to show that for any $\mathscr{C}$-quasi-isomorphism $s: Y^\bullet\to X^\bullet$ with
$X^\bullet\in K^{-}(\mathscr{A})$, there exists a morphism $f: Z^\bullet\to Y^\bullet$ with
$Z^\bullet\in K^{-}(\mathscr{A})$ such that $sf$ is a $\mathscr{C}$-quasi-isomorphism.

Suppose $X^n\neq 0$ with $X^i=0$ for any $i>n$. Then there exists a commutative diagram:
$$\xymatrix{Z^\bullet : \ar[d]^{f} & \cdots\ar[r] & Y^{n-1}\ar[r]\ar@{=}[d] & Y^{n}\ar[r]\ar@{=}[d] & \Ker d^{n+1}_Y\ar[d]\ar[r] & 0\\
Y^\bullet : \ar[d]^{s} & \cdots\ar[r] & Y^{n-1}\ar[r]\ar[d] & Y^{n}\ar[r]\ar[d] & Y^{n+1}\ar[d]\ar[r] & \cdots\\
X^\bullet :  & \cdots\ar[r] & X^{n-1}\ar[r] & X^{n}\ar[r] & 0\ar[r] & \cdots,}
$$
where $\Ker d^{n+1}_Y\to Y^{n+1}$ is the canonical map. Since both $f$ and $s$ are $\mathscr{C}$-quasi-isomorphisms, so is $sf$.

(2) Consider the canonical map $G:\Hom_{K(\mathscr{A})}(C^{\bullet},X^{\bullet})\to \Hom_{D_\mathscr{C}(\mathscr{A})}(C^{\bullet},X^{\bullet})$
defined by $G(f)=f/\id_{C^\bullet}$. If $G(f)=0$, then there exists a $\mathscr{C}$-quasi-isomorphism $s:Z^\bullet\to C^\bullet$ such that
$fs\sim 0$. By Lemma 2.4(1) there exists a cochain map $g:C^\bullet\to Z^\bullet $ such that $sg\sim \id_{C^\bullet}$, and then $f\sim 0$.
On the other hand, let $f/s\in \Hom_{D_\mathscr{C}(\mathscr{A})}(C^{\bullet},X^{\bullet})$, that is, it has a diagram of the form
$C^\bullet\buildrel {s} \over \Longleftarrow Z^\bullet\buildrel {f} \over\longrightarrow X^\bullet$, where $s$ is a $\mathscr{C}$-quasi-isomorphism.
It follows from Lemma 2.4(1) there exists a cochain map $g:C^\bullet\to Z^\bullet $ such that $sg\sim \id_{C^\bullet}$,
which implies that $f/s=(fg)/\id_{C^\bullet}=G(fg)$. Thus $G$ is an isomorphism, as desired.

(3) Let $F:\mathscr{A}\to  D_{\mathscr{C}}^{b}(\mathscr{A})$ denote the composition functor, it suffices to show that for any $M,N\in \mathscr{A}$,
the map $F:\Hom_\mathscr{A}(M,N)\to \Hom_{ D_{\mathscr{C}}^{b}(\mathscr{A})}(M,N)$ is an isomorphism.

Let $f\in \Hom_\mathscr{A}(M,N)$. If $F(f)=0$, then there exists a $\mathscr{C}$-quasi-isomorphism $s:Z^\bullet\to M$ such that $fs\sim 0$,
and then $H^0(f)H^0(s)=0$. Since $H^0(s)$ is an isomorphism, $f=H^0(f)=0$. On the other hand, let $f/s\in \Hom_{D_\mathscr{C}^b(\mathscr{A})}(M,N)$,
that is, it has a diagram of the form $M\buildrel {s} \over \Longleftarrow Z^\bullet\buildrel {f} \over\longrightarrow N$, where $s$ is a
$\mathscr{C}$-quasi-isomorphism. Then $H^0(s):H^0(Z^\bullet)\to M$ is an isomorphism. Put $g:=H^0(f)H^0(s)^{{\text-}1}\in \Hom_\mathscr{A}(M,N)$.
Consider the truncation:
$$U^\bullet:= \cdots\to Z^{-2}\buildrel {d_Z^{-2}} \over\longrightarrow Z^{-1}\buildrel {d_Z^{-1}} \over\longrightarrow \Ker d^0\to 0 $$
of $Z^\bullet$ and the canonical map $i:U^\bullet\to Z^\bullet$. Since $s$ is a $\mathscr{C}$-quasi-isomorphism,
so is $si$. We have the following commutative diagram:
$$\xymatrix{U^\bullet\ar[r]^{i}\ar[d] & Z^\bullet\ar[d]^{s}\\
H^0(Z^\bullet)\ar[r]^{H^0(s)} & M,}
$$
where $U^\bullet\to H^0(Z^\bullet)$ is the canonical map,
so $gsi=H^0(f)H^0(s)^{{\text-}1}si=fi$. Then we get the following commutative diagram of complexes:
$$\xymatrix{&Z^\bullet\ar@{=>}[dl]_s\ar[dr]^f&\\
M & U^\bullet\ar@{=>}[l]^{si}\ar[u]^i\ar[d]^{si}\ar[r]^{fi} & N\\
 & M,\ar@{=>}[ul]^{\id_M}\ar[ur]_g & }$$
which implies $F(g)=g/\id_M=f/s$.
\hfill{$\square$}

\vspace{0.2cm}

Set $K^{-,\mathscr{C}b}(\mathscr{C}):=\{X^{\bullet}\in K^{-}(\mathscr{C})\mid$ there exists $n\in \mathbb{Z}$ such that
$H^i(\Hom_{\mathscr{A}}(C,X^{\bullet}))=0$ for any $C\in\mathscr{C}$ and $i \leq n\}$.

\vspace{0.2cm}

{\bf Proposition 3.4.} {\it ([AHV, Theorem 3.3]) If $\mathscr{C}$ is a contravariantly finite subcategory of $\mathscr{A}$,
then we have a triangle-equivalence $K^{-,\mathscr{C}b}(\mathscr{C})\cong D_{\mathscr{C}}^{b}(\mathscr{A})$.}

\vspace{0.2cm}
In the rest of this section, we always suppose that $\mathscr{C}$ is a contravariantly finite subcategory of $\mathscr{A}$ unless otherwise specified.

\vspace{0.2cm}
{\bf Definition 3.5.} Let $\mathscr{D}$ be a subclass of objects in $\mathscr{A}$ and $M\in \mathscr{A}$.

(1) A {\it $\mathscr{C}$-proper $\mathscr{D}$-resolution} of $M$ is a $\mathscr{C}$-quasi-isomorphism $f: D^{\bullet}\to M$,
where $D^{\bullet}$ is a complex of objects in $\mathscr{D}$ with $D^n=0$ for any $n>0$, that is, it has an associated
$\Hom_\mathscr{A}(\mathscr{C},-)$-exact complex
$\cdots \to D^{-n}\to D^{-n+1}\to \cdots\to D^0\buildrel {f} \over \longrightarrow M\to 0$.

(2) The {\it $\mathscr{C}$-proper $\mathscr{D}$-dimension} of $M$, written $\mathscr{C}\mathscr{D}$-$\dim M$,
is defined as $\inf\{n\mid$ there exists a $\Hom_\mathscr{A}(\mathscr{C},-)$-exact complex
$0 \to D^{-n}\to D^{-n+1}\to \cdots\to D^0\buildrel {f} \over \longrightarrow M\to 0$\}.
If no such an integer exists, then set $\mathscr{C}\mathscr{D}$-$\dim M=\infty$.

(3) For a class $\mathscr{E}$ of objects of $\mathscr{A}$, the {\it $\mathscr{C}$-proper $\mathscr{D}$-dimension} of $\mathscr{E}$,
written $\mathscr{C}\mathscr{D}$-$\dim \mathscr{E}$, is defined as $\sup\{\mathscr{C}\mathscr{D}$-$\dim M\mid M\in \mathscr{E}\}$.


\vspace{0.2cm}
{\bf Remark 3.6.} (1) If $\mathscr{A}$ has enough projective objects and $\mathscr{C}=\mathscr{P}(\mathscr{A})$,
then a $\mathscr{C}$-proper $\mathscr{D}$-resolution is just a $\mathscr{D}$-resolution and the $\mathscr{C}$-proper
$\mathscr{D}$-dimension of an object $M\in \mathscr{A}$ is just the usual $\mathscr{D}$-dimension $\mathscr{D}$-$\dim M$ of $M$.

(2) If $\mathscr{D}=\mathscr{C}$, then a $\mathscr{C}$-proper $\mathscr{D}$-resolution is just a $\mathscr{C}$-proper resolution.
In this case, it is also called a {\it left $\mathscr{C}$-resolution} and the $\mathscr{C}$-proper $\mathscr{D}$-dimension is
the left $\mathscr{C}$-dimension (see [EJ2]).

\vspace{0.2cm}

Let $M \in\mathscr{A}$. Since $\mathscr{C}$ is a contravariantly finite subcategory of $\mathscr{A}$, we may choose
a left $\mathscr{C}$-resolution $C^{\bullet}_M\to M$ of $M$. Put $\Ext_{\mathscr{C}}^n(M,N):=H^{n}\Hom_\mathscr{A}(C^{\bullet}_M,N)$
for any $N\in \mathscr{A}$ and $n \in \mathbb{Z}$.
Note that
$C^{\bullet}_M$ is isomorphic to $M$ in $D_\mathscr{C}(\mathscr{A})$. By Proposition 3.3(1)(2), we have
$\Ext_{\mathscr{C}}^n(M,N)=H^{n}\Hom_\mathscr{A}(C^{\bullet}_M,N)=\Hom_{K(\mathscr{A})}(C^{\bullet}_M,N[n])\cong
\Hom_{D_\mathscr{C}(\mathscr{A})}(C^{\bullet}_M,N[n])\cong\Hom_{D_\mathscr{C}^b(\mathscr{A})}(M,N[n])$.

\vspace{0.2cm}

The following is cited from [EJ2, Chapter 8].

\vspace{0.2cm}

{\bf Lemma 3.7.} (1) For any $M\in \mathscr{A}$, the functor $\Ext_{\mathscr{C}}^n(M,-)$ does not
depend on the choices of left $\mathscr{C}$-resolutions of $M$.

(2) For any $M\in \mathscr{A}$ and $n<0$, $\Ext_{\mathscr{C}}^n(M,-)=0$ and there exists a natural equivalence
$\Hom_\mathscr{A}(M,-) \cong \Ext_{\mathscr{C}}^0(M,-)$ whenever $\mathscr{C}$ is admissible.

(3) If $\mathscr{C}$ is admissible, then every $\Hom_\mathscr{A}(\mathscr{C},-)$-exact complex $0\to L\to M\to N\to 0$ induces a long exact sequence
$0\to \Hom_\mathscr{A}(N,-)\to \Hom_\mathscr{A}(M,-)\to \Hom_\mathscr{A}(L,-)\to \cdots \to \Ext_{\mathscr{C}}^n(N,-)\to
\Ext_{\mathscr{C}}^n(M,-)\to \Ext_{\mathscr{C}}^n(L,-)\to \Ext_{\mathscr{C}}^{n+1}(N,-)\to \cdots$.

\vspace{0.2cm}
{\bf Theorem 3.8.} {\it Let $\mathscr{C}$ be admissible and closed under direct summands, then the following statements are equivalent
for any $M \in\mathscr{A}$ and $n\geq 0$.

(1) $\mathscr{C}\mathscr{C}$-$\dim M \leq n$.

(2) $\Ext_{\mathscr{C}}^i(M,N)=0$ for any $N\in \mathscr{A}$ and $i\geq n+1$.

(3) $\Ext_{\mathscr{C}}^{n+1}(M,N)=0$ for any $N\in \mathscr{A}$.

(4) For any left $\mathscr{C}$-resolution $C^{\bullet}_M\to M$ of $M$, $\Ker d_{C_M}^{-n+1}\in\mathscr{C}$,
where $d_{C_M}^{-n+1}$ is the $(-n+1)$st differential of $C^{\bullet}_M$.}

\vspace{0.2cm}

{\it Proof.} $(1)\Rightarrow (2)$ Let $0 \to C^{-n}\to C^{-n+1}\to \cdots\to C^0\to M\to 0$ be a
left $\mathscr{C}$-resolution of $M$. Then $\Hom_\mathscr{A}(C^{-i},N)$ =0 for any $N\in \mathscr{A}$ and $i\geq n+1$
and the assertion follows.

$(2)\Rightarrow (3)$ and $(4)\Rightarrow (1)$ are trivial.

$(3)\Rightarrow (4)$ Let $\cdots \to C^{-n}_M\buildrel {d_{C_M}^{-n}} \over \longrightarrow C^{-n+1}_M\to \cdots\to C^0_M\to M\to 0$
be a left $\mathscr{C}$-resolution of $M$. Then we get a $\Hom_\mathscr{A}(\mathscr{C},-)$-exact exact sequence
$0\to \Ker d_{C_M}^{-n}\to C^{-n}_M\to \Ker d_{C_M}^{-n+1}\to 0$. Since $\Ext_{\mathscr{C}}^{n+1}(M,\Ker d_{C_M}^{-n})=0$,
$\Ext_{\mathscr{C}}^{1}(\Ker d_{C_M}^{-n+1}, \Ker d_{C_M}^{-n})\cong \Ext_{\mathscr{C}}^{n+1}(M,\Ker d_{C_M}^{-n})=0$ by the
dimension shifting. Applying $\Hom_\mathscr{A}(-,\Ker d_{C_M}^{-n})$ to the exact sequence
$0\to \Ker d_{C_M}^{-n}\to C^{-n}_M\to \Ker d_{C_M}^{-n+1}\to 0$, it follows from Lemma 3.7(3) that the sequence splits.
So $\Ker d_{C_M}^{-n+1}$ is a direct summand of $C^{-n}_M$ and $\Ker d_{C_M}^{-n+1}\in\mathscr{C}$. \hfill{$\square$}

\vspace{0.5cm}

\centerline{\bf 4. $\mathscr{C}$-singularity categories}

\vspace{0.2cm}

In this section, unless otherwise specified, we always suppose that $A$ is a finite-dimensional algebra over a field, $\mathscr{A}=A{\text-}\mod$
and $\mathscr{C}$ is a full and additive subcategory of $\mathscr{A}$ which is contravariantly finite in $\mathscr{A}$ and is admissible and closed under direct summands.

Recall that an additive category is called a {\it Krull-Schmidt category} if each of its object $X$ has a decomposition $X\cong X_1\bigoplus X_2\bigoplus \cdots \bigoplus X_n$
such that each $X_i$ is indecomposable with a local endomorphism ring. By [BD, Proposition A.2] $K^{b}(\mathscr{C})$ is a Krull-Schmidt category,
so it is closed under direct summands and $K^{b}(\mathscr{C})$ viewed as a full triangulated subcategory of $D^{b}_\mathscr{C}(\mathscr{A})$ is thick.
It is of interest to consider the quotient triangulated category $D^{b}_\mathscr{C}(\mathscr{A})$ /$K^{b}(\mathscr{C})$.

\vspace{0.2cm}

{\bf Definition 4.1.} We call $D_{\mathscr{C}{\text-}sg}(\mathscr{A}):= D^{b}_\mathscr{C}(\mathscr{A})$ /$K^{b}(\mathscr{C})$ the {\it $\mathscr{C}$-singularity category}.

\vspace{0.2cm}

{\bf Example 4.2.} (1) If $\mathscr{C}=A$-$\proj$, then $D^{b}_\mathscr{C}(\mathscr{A})$ is the usual bounded derived category $D^{b}(\mathscr{A})$ and
the $\mathscr{C}$-singularity category $D_{\mathscr{C}{\text-}sg}(\mathscr{A})$ is the singularity category $D_{sg}(A)$ which is called the ``stabilized derived category" in [Bu].

(2) Let $\mathscr{C}=\G(A)$ (the subcategory of $A$-$\mod$ consisting of Gorenstein projective modules). If $\G(A)$ is contravariantly finite in $A$-$\mod$,
for example, if $A$ is Gorenstein (that is, the left and right self-injective dimensions of $A$ are finite) or $\G(A)$ contains only finitely many non-isomorphic indecomposable modules,
then the bounded $\mathscr{C}$-derived category of $\mathscr{A}$, denoted by $D^{b}_{\mathscr{G}(A)}(\mathscr{A})$,
is the {\it bounded Gorenstein derived category} introduced in [GZ]. The $\mathscr{C}$-singularity category $D_{\mathscr{G}(A){\text-}sg}(\mathscr{A})$
is the quotient triangulated category $D^{b}_{\mathscr{G}(A)}(\mathscr{A})$ /$K^{b}(\mathscr{G}(A))$, we call it the {\it Gorenstein singularity category}.

\vspace{0.2cm}

Given a complex $X^\bullet$ and an integer $i\in \mathbb{Z}$, we denote by $\sigma^{\geq i}X^\bullet$ the complex with $X^j$ in the $j$th degree
whenever $j\geq i$ and 0 elsewhere, and set $\sigma^{>i}X^\bullet:=\sigma^{\geq i+1}X^\bullet$. Dually, for the notations $\sigma^{\leq i}X^\bullet$ and $\sigma^{< i}X^\bullet$.
Recall that the cardinal of the set $\{X^i \neq 0\mid i\in \mathbb{Z}\}$ is called the {\it width} of $X^\bullet$, and denoted by $\omega(X^\bullet)$.

It is well known that $A$ has finite global dimension if and only if $D_{sg}(A)=0$. For the $\mathscr{C}$-singularity category $D^{b}_{\mathscr{C}{\text-}sg}(\mathscr{A})$
we have the following property.

\vspace{0.2cm}
{\bf Proposition 4.3.} {\it If $\mathscr{C}\mathscr{C}$-dim $\mathscr{A}< \infty$, then $D_{\mathscr{C}{\text-}sg}(\mathscr{A})=0$.}

\vspace{0.2cm}
{\it Proof.} We claim that for every $X^\bullet\in K^{b}(\mathscr{A})$ there exists a $\mathscr{C}$-quasi-isomorphism $C_X^\bullet \to X^\bullet$ such that
$C_X^\bullet \in K^{b}(\mathscr{C})$. We proceed
by induction on the width $\omega(X^\bullet)$ of $X^\bullet$.

Let $\omega(X^\bullet)$=1. Because $\mathscr{C}$ is contravariantly finite and $\mathscr{C}\mathscr{C}$-$\dim \mathscr{A}< \infty$,
there exists a $\mathscr{C}$-quasi-isomorphism $C_X^\bullet \to X^\bullet$ with $C_X^\bullet \in K^{b}(\mathscr{C})$.

Let $\omega(X^\bullet)\geq 2$ with $X^j\neq 0$ and $X^i=0$ for any $i<j$. Put $X^\bullet_1$:=$X^j[-j-1]$, $X^\bullet_2:=\sigma^{> j}X^\bullet$ and $g=d_X^j[-j-1]$.
We have a distinguished triangle $X^\bullet_1\buildrel {g} \over \longrightarrow X^\bullet_2\to X^\bullet\to X^\bullet_1[1]$ in $K^b(\mathscr{A})$.
By the induction hypothesis, there exist $\mathscr{C}$-quasi-isomorphisms
$f_{X_1}$: $C_{X_1}^\bullet\to X^\bullet_1$ and $f_{X_2}$: $C_{X_2}^\bullet\to X^\bullet_2$ with $C_{X_1}^\bullet$,$C_{X_2}^\bullet$ $\in K^{b}(\mathscr{C})$.
Then by Remark 2.2(1) and Lemma 2.3, $f_{X_2}$ induces an isomorphism:
$$\Hom_{K^b(\mathscr{A})}(C_{X_1}^\bullet,C_{X_2}^\bullet)\cong \Hom_{K^b(\mathscr{A})}(C_{X_1}^\bullet,X^\bullet_2).$$
So there exists a morphism $f^\bullet:C_{X_1}^\bullet \to C_{X_2}^\bullet$, which is unique up to homotopy, such that
$f_{X_2} f^\bullet=gf_{X_1}$. Put $C_{X}^\bullet=\Con(f^\bullet)$. We have the following distinguished triangle in $K^b(\mathscr{C})$:
\begin{center}
$C_{X_1}^\bullet \buildrel {f^\bullet} \over \longrightarrow C_{X_2}^\bullet\to C_{X}^\bullet\to C_{X_1}^\bullet[1].$
\end{center}
Then there exists a morphism $f_{X}:C_{X}^\bullet\to X^\bullet$ such that the following diagram commutes:
$$\xymatrix{C_{X_1}^\bullet \ar[r]^{f^\bullet}\ar[d]^{f_{X_1}} & C_{X_2}^\bullet \ar[r] \ar[d]^{f_{X_2}} &
C_{X}^\bullet \ar[r]\ar@{-->}[d]^{f_{X}} & C_{X_1}^\bullet[1] \ar[d]^{f_{X_1}[1]} \\
X^\bullet_1 \ar[r]^{g} & X^\bullet_2 \ar[r] &
X^\bullet \ar[r] & {X^\bullet_1}[1].}$$
For any $C\in \mathscr{C}$ and any $n\in\mathbb{Z}$, we have the following commutative diagram with exact rows:
$$\xymatrix{(C,C_{X_1}^\bullet[n]) \ar[r]\ar[d]^{(C, f_{X_1}[n])} & (C,C_{X_2}^\bullet[n]) \ar[r] \ar[d]^{(C, f_{X_2}[n])} &
(C,C_{X}^\bullet[n]) \ar[r]\ar@{-->}[d]^{(C, f_{X}[n])} & (C,C_{X_1}^\bullet[n+1]) \ar[d]^{(C, f_{X_1}[n+1])}\ar[r] & (C,C_{X_2}^\bullet[n+1]) \ar[d]^{(C,f_{X_2}[n+1])}\\
(C,X^\bullet_1[n]) \ar[r] & (C,X^\bullet_2[n]) \ar[r] &
(C,X^\bullet[n]) \ar[r] & (C,{X^\bullet_1}[n+1])\ar[r] & (C,{X^\bullet_2}[n+1]),}$$
where $(C,-)$ denotes the functor $\Hom_{K(\mathscr{A})}(C,-)$. Since $f_{X_1}$ and $f_{X_2}$ are $\mathscr{C}$-quasi-isomorphisms,
$(C,f_{X_1}[n])$ and $(C,f_{X_2}[n])$ are isomorphisms, and hence so is $(C,f_{X}[n])$ for each $n$, that is, $f_{X}$ is a $\mathscr{C}$-quasi-isomorphism.
The claim is proved.

It follows from the claim that every object $X^\bullet$ in $D^{b}_\mathscr{C}(\mathscr{A})$ is isomorphic to some $C_X^\bullet$ of $K^{b}(\mathscr{C})$
in $D^{b}_\mathscr{C}(\mathscr{A})$. Thus $D_{\mathscr{C}{\text-}sg}(\mathscr{A})=0$. \hfill{$\square$}

\vspace{0.2cm}

As an application of Proposition 4.3, we have the following

\vspace{0.2cm}

{\bf Corollary 4.4.} {\it
(1) $\mathscr{C}\mathscr{C}$-$\dim M< \infty$ for any $M\in \mathscr{A}$ if and only if
$D_{\mathscr{C}{\text-}sg}(\mathscr{A})=0$.

(2) If $A$ is of finite representation type, then $\mathscr{C}\mathscr{C}$-$\dim \mathscr{A}< \infty$ if and only if
$D_{\mathscr{C}{\text-}sg}(\mathscr{A})=0$.}

\vspace{0.2cm}

{\it Proof.} In both assertions, the necessity follows from Proposition 4.3. In the following, we only need to prove the sufficiency.

(1) Let $D_{\mathscr{C}{\text-}sg}(\mathscr{A})=0$ and $M\in \mathscr{A}$.
Then $M=0$ in $D_{\mathscr{C}{\text-}sg}(\mathscr{A})$ and
$M$ is isomorphic to $C^\bullet$ in $D_{\mathscr{C}}^b(\mathscr{A})$
for some $C^\bullet\in K^b(\mathscr{C})$. We use the equivalent class of right fractions to denote a morphism
in $D_{\mathscr{C}}^b(\mathscr{A})$. Let $f/s: C^\bullet \buildrel {s} \over \Longleftarrow Z^\bullet \buildrel {f} \over \longrightarrow M$
be an isomorphism in $D_{\mathscr{C}}^b(\mathscr{A})$, where $s$ is a $\mathscr{C}$-quasi-isomorphism. Then $f$ is a $\mathscr{C}$-quasi-isomorphism.
By Lemma 2.4(1), there exists a $\mathscr{C}$-quasi-isomorphism $s^{'}:C^\bullet\to Z^\bullet$. So $fs^{'}:C^\bullet\to M$ is also a
$\mathscr{C}$-quasi-isomorphism
and hence $H^i\Hom_{\mathscr{A}}(C,C^\bullet)=0$ whenever $C\in \mathscr{C}$ and $i\neq 0$. Consider the truncation:
$${C^{'}}^\bullet:= \cdots\to C^{-2}\to C^{-1}\to \Ker d_C^0\to 0$$
of $C^\bullet$. Then the composition ${C^{'}}^\bullet\hookrightarrow C^\bullet\buildrel {fs^{'}}\over\longrightarrow M$ is a $\mathscr{C}$-quasi-isomorphism.
Notice that $C^\bullet\in K^b(\mathscr{C})$, we may suppose $C^n\neq 0$ and $C^i=0$ whenever $i>n$. Then we have a $\mathscr{C}$-acyclic complex
$0\to \Ker d_C^0\to C^0\buildrel {d_C^0}\over \longrightarrow C^1 \to \cdots \to C^n \to 0$ with all $C^i$ in $\mathscr{C}$.
Because $\mathscr{C}$ is closed under direct summands, $\Ker d_C^0\in\mathscr{C}$ and $\mathscr{C}\mathscr{C}$-$\dim M< \infty$.

(2) Let $A$ be of finite representation type, and let $\{M_i\mid 1\leq i\leq n\}$ be the set of all non-isomorphic
indecomposable modules in $\mathscr{A}$. By (1) $\mathscr{C}\mathscr{C}$-$\dim M_i< \infty$ for any $1\leq i\leq n$. Now set
$m=\sup\{\mathscr{C}\mathscr{C}$-$\dim M_i\mid 1\leq i\leq n\}$. Since $\mathscr{A}$ is Krull-Schmidt, every module $M\in \mathscr{A}$
can be decomposed into a finite direct sum of modules in $\{M_i\mid 1\leq i\leq n\}$. Then it is easy to see that
$\mathscr{C}\mathscr{C}$-$\dim M\leq m$ and $\mathscr{C}\mathscr{C}$-$\dim\mathscr{A}\leq m< \infty$.
\hfill{$\square$}

\vspace{0.2cm}

As a consequence of Corollary 4.4(1), we have the following


\vspace{0.2cm}

{\bf Corollary 4.5.} {\it If $A$ is Gorenstein, then $D_{\mathscr{G}(A){\text-}sg}(\mathscr{A})=0$.}

\vspace{0.2cm}

{\it Proof.} Let $A$ be Gorenstein. Because $A$-$\proj\subseteq\mathscr{G}(A)$, we have that $\mathscr{G}(A)$
is admissible in $A$-$\mod$ by [EJ2, Remark 11.5.2]. By [Hos, Theorem], we have $\G(A)$-$\dim M<\infty$ for any
$M\in \mathscr{A}$. So $D_{\mathscr{G}(A){\text-}sg}(\mathscr{A})=0$ by [AvM, Proposition 4.8] and Corollary 4.4(1).
\hfill{$\square$}

\vspace{0.2cm}


Put $\mathscr{G}(\mathscr{C})=\{M\cong \Im(C^{-1}\to C^0)\mid$ there exists an acyclic complex $\cdots\to C^{-1}\to C^0\to C^1\to
\cdots$ in $\mathscr{C},\text{ which is both }\Hom_\mathscr{A}(\mathscr{C},-)\text{-exact and }\Hom_\mathscr{A}(-,\mathscr{C})\text{-exact}\}$,
see [SSW], where it is called the {\it Gorenstein category} of $\mathscr{C}$. This notion unifies the following ones:
modules of Gorenstein dimension zero ([AuB]), Gorenstein projective
modules, Gorenstein injective modules ([EJ1]), $V$-Gorenstein
projective modules, $V$-Gorenstein injective modules ([EJL]), and so on.
Set $\mathscr{G}^1(\mathscr{C})=\mathscr{G}(\mathscr{C})$ and
inductively set $\mathscr{G}^n(\mathscr{C})=\mathscr{G}(\mathscr{G}^{n-1}(\mathscr{C}))$ for any $n\geq 2$.
It was shown in [SSW] that $\mathscr{G(C)}$ possesses many nice properties when $\mathscr{C}$ is self-orthogonal. For example, in this case,
$\mathscr{G(C)}$ is closed under extensions and $\mathscr{C}$ is a projective generator and an injective cogenerator for $\mathscr{G(C)}$,
which induce that $\mathscr{G}^n(\mathscr{C})=\mathscr{G(C)}$ for any $n \geq 1$, see [SSW] for more details. Later on, Huang
generalized this result to an arbitrary full and additive subcategory $\mathscr{C}$ of $\mathscr{A}$, see [Hu].



Denote by $\varepsilon$ the class of all $\Hom_\mathscr{A}(\mathscr{C},-)$-exact complexes of the form: $0\to L\buildrel {i} \over
\longrightarrow M\buildrel {p} \over \longrightarrow N\to 0$ with $L,M,N\in \mathscr{G}(\mathscr{C})$. We have the following fact.

\vspace{0.2cm}

{\bf Proposition 4.6.} {\it $(\mathscr{G(C)},\varepsilon)$  is an exact category.}

\vspace{0.2cm}
{\it Proof.} We will prove that all the axioms in Definition 2.6 are satisfied.
It is trivial that the axiom [E0] is satisfied. In the following, we prove that the other axioms are satisfied.

For [${\text E1^{op}}$], let  $f:G_1\to G_2$ and  $g:G_2\to G_3$ be admissible epics in $\mathscr{G(C)}$. Then it is easy to see that $gf$
is also an admissible epic. By Lemma 3.7(3), the following $\Hom_\mathscr{A}(\mathscr{C},-)$-exact sequence:
$$0\to \Ker gf\to G_1\buildrel {gf} \over \longrightarrow G_3\to 0$$
is also $\Hom_\mathscr{A}(-,\mathscr{C})$-exact. It follows from [Hu, Proposition 4.7] that $\Ker gf\in \mathscr{G(C)}$.

For [${\text E2^{op}}$], let $f:G_2\to G_3$ be an admissible epic in $\mathscr{G(C)}$ and $g:G^{'}_2\to G_3$ an arbitrary morphism in $\mathscr{G(C)}$. We have the following pull-back diagram
with the second row in $\varepsilon$:
$$\xymatrix{0\ar[r]&G_1\ar[r]^{h^{'}}\ar@{=}[d]&X\ar[r]^{f^{'}}\ar[d]^{g^{'}}&G^{'}_2\ar[d]^g\ar[r]&0\\
0\ar[r]&G_1\ar[r]^{h}&G_2\ar[r]^{f}&G_3\ar[r]&0.
}$$
For any $C\in \mathscr{C}$ and any morphism $\varphi: C\to G^{'}_2$, there exists a morphism $\phi: C\to G_2$ such that $g\varphi=f\phi$. Notice that the right square is a pull-back diagram,
so there exists a morphism $\phi^{'}: C\to X$ such that $\varphi=f^{'}\phi^{'}$
and hence the exact sequence $0\to G_1\buildrel {h^{'}} \over \longrightarrow X\buildrel {f^{'}} \over \longrightarrow G^{'}_2\to 0$ is
$\Hom_\mathscr{A}(\mathscr{C},-)$-exact. It follows from Lemma 3.7(3) that this sequence is also $\Hom_\mathscr{A}(-,\mathscr{C})$-exact.
By [Hu, Proposition 4.7], $X\in \mathscr{G(C)}$, which implies that $0\to G_1\buildrel {h^{'}} \over \longrightarrow X\buildrel {f^{'}} \over \longrightarrow G^{'}_2\to 0$
lies in $\varepsilon$.

For [E2], let $f:G_1\to G_2$ be an admissible monic in $\mathscr{G(C)}$ and $g:G_1\to G^{'}_2$ an arbitrary morphism in $\mathscr{G(C)}$. We have the following push-out diagram
with the first row in $\varepsilon$:
$$\xymatrix{0\ar[r]&G_1\ar[r]^f\ar[d]^g&G_2\ar[r]^h\ar[d]^{g^{'}}&G_3\ar@{=}[d]\ar[r]&0\\
0\ar[r]&G^{'}_2\ar[r]^{f^{'}}&D\ar[r]^{h^{'}}&G_3\ar[r]&0.
}$$
For any $C\in \mathscr{C}$ and any morphism $\varphi: C\to G_3$, there exists a morphism $\phi: C\to G_2$ such that $\varphi=h\phi=h^{'}g^{'}\phi$.
So the exact sequence $0\to G^{'}_2\buildrel {f^{'}} \over \longrightarrow D\buildrel {h^{'}} \over \longrightarrow G_3\to 0$ is
$\Hom_\mathscr{A}(\mathscr{C},-)$-exact. It follows from Lemma 3.7(3) that this sequence is also $\Hom_\mathscr{A}(-,\mathscr{C})$-exact.
By [Hu, Proposition 4.7], $D\in \mathscr{G(C)}$, which implies that $0\to G^{'}_2\buildrel {f^{'}} \over \longrightarrow D\buildrel {h^{'}} \over \longrightarrow G_3\to 0$
lies in $\varepsilon$ .

Now let $0\to G_0\buildrel {i} \over \longrightarrow G_1\to G_2\to0$ and $0\to G_1\buildrel {j} \over \longrightarrow  G^{'}_1\to G^{''}_1\to 0$ lie in $\varepsilon$.
We have the following push-out diagram:

$$\xymatrix{ & & 0\ar[d]& 0\ar[d]& \\
0\ar[r] & G_0\ar@{=}[d]\ar[r]^i & G_1\ar[d]^j\ar[r]& G_2 \ar[d]\ar[r]& 0\\
0\ar[r] & G_0\ar[r]^{ji} & G^{'}_1\ar[d]\ar[r]& G^{'}_2 \ar[d]\ar[r]& 0\\
 & & G^{''}_1\ar[d]\ar@{=}[r] & G^{''}_1\ar[d]& \\
 & & 0 & 0.& }$$
By [E2], the rightmost column lies in $\varepsilon$. For any $C\in \mathscr{C}$, applying the functor $(C,-):=\Hom_\mathscr{A}(C,-)$ to the commutative diagram
we get the following commutative diagram:
$$\xymatrix{ & & 0\ar[d]& 0\ar[d]& \\
0\ar[r] & (C,G_0)\ar@{=}[d]\ar[r]^{(C,i)} & (C,G_1)\ar[d]^{(C,j)}\ar[r]& (C,G_2) \ar[d]\ar[r]& 0\\
0\ar[r] & (C,G_0)\ar[r]^{(C,ji)} & (C,G^{'}_1)\ar[d]\ar[r]& (C,G^{'}_2) \ar[d]& \\
 & & (C,G^{''}_1)\ar[d]\ar@{=}[r] & (C,G^{''}_1)\ar[d]& \\
 & & 0 & 0. & }$$
By the snake lemma, the morphism $(C,G^{'}_1)\to (C,G^{'}_2)$ is epic. Then $0\to G_0\buildrel {ji} \over \longrightarrow G^{'}_1\to G^{'}_2 \to 0$ lies
in $\varepsilon$, and [E1] follows. \hfill{$\square$}

\vspace{0.2cm}

By Proposition 4.6, we have the following

\vspace{0.2cm}

{\bf Corollary 4.7.} $(\mathscr{G(C)},\varepsilon)$ is a Frobenius category, that is, $(\mathscr{G(C)},\varepsilon)$ has enough projective objects and enough
injective objects such that the projective objects coincide with the injective objects.

\vspace{0.2cm}

{\it Proof.}  Because $\mathscr{C}$ is the class of (relative) projective-injective
objects in $\mathscr{G(C)}$, the assertion follows from Proposition 4.6. \hfill{$\square$}

\vspace{0.2cm}

For $M, N\in \mathscr{A}$, let $\mathscr{C}(M,N)$ denote the subspace of $A$-maps from $M$ to $N$ factoring through $\mathscr{C}$.
Put $^{\bot_\mathscr{C}}\mathscr{C}=\{M\in \mathscr{A}\mid \Ext_\mathscr{C}^i(M,C)=0$ for any $C\in\mathscr{C}$ and $i\geq1\}$. By definition,
it is clear that $\mathscr{C}\subseteq\mathscr{G(C)}\subseteq {^{\bot_\mathscr{C}}\mathscr{C}}$.

\vspace{0.2cm}

{\bf Lemma 4.8.} {\it For any $M\in {^{\bot_\mathscr{C}}\mathscr{C}}$ and $N\in \mathscr{A}$, we have  a canonical isomorphism of abelian groups:
$$\Hom_\mathscr{A}(M,N)/\mathscr{C}(M,N)\cong \Hom_{D_{\mathscr{C}{\text-}sg}(\mathscr{A})}(M,N).$$}

\vspace{0.2cm}

{\it Proof.} In the following, a morphism from $M$ to $N$ in $D_{\mathscr{C}{\text-}sg}(\mathscr{A})$ is denoted by
the equivalent class of left fractions $s\backslash a: M\buildrel {a} \over \longrightarrow Z^\bullet \buildrel {s} \over \Longleftarrow N$, where
$Z^\bullet\in D^b_\mathscr{C}(\mathscr{A})$ and $\Con(s)\in K^b(\mathscr{C})$. We have a distinguished triangle in $D^b_\mathscr{C}(\mathscr{A})$:
$$N\buildrel {s} \over \Longrightarrow Z^\bullet\to \Con(s)\to N[1].\eqno{(1)}$$
Consider the canonical map $G : \Hom_\mathscr{A}(M,N)\to \Hom_{D_{\mathscr{C}{\text -}sg}(\mathscr{A})}(M,N)$ defined by $G(f)=\id_N\backslash f$.
We first prove that $G$ is surjective. For any $N\in \mathscr{A}$, we have the following left $\mathscr{C}$-resolution of $N$:
$$\cdots\to C^{-n}\buildrel {d_C^{-n}} \over \longrightarrow C^{-n+1}\to \cdots\buildrel {d_C^{-1}} \over \longrightarrow C^0 \buildrel {d_C^0} \over \longrightarrow N\to 0.$$
Then in $D_\mathscr{C}(\mathscr{A})$, $N$ is isomorphic to  the complex $C^\bullet:= \cdots\to C^{-n}\buildrel {d_C^{-n}} \over \longrightarrow C^{-n+1}\to
\cdots\buildrel {d_C^{-1}} \over \longrightarrow C^0 \to 0$, and so is isomorphic to the complex $0\to \Ker d_C^{-l}\to C^{-l}\buildrel {d_C^{-l}} \over
\longrightarrow C^{-l+1}\to \cdots\buildrel {d_C^{-1}} \over \longrightarrow C^0 \to 0$ for any $l\geq 0$. Hence we have a distinguished triangle in $D^b_\mathscr{C}(\mathscr{A})$:
$$\Ker d_C^{-l}[l]\to \sigma^{\geq -l}C^\bullet\buildrel {d_C^0} \over \longrightarrow N\buildrel {s'} \over \Longrightarrow \Ker d_C^{-l}[l+1],\eqno{(2)}$$
where $\Con(s')\in K^b(\mathscr{C})$. Since  $\Con(s)\in K^b(\mathscr{C})$, it follows from Proposition 3.3 that there exists $l_0\gg 0$ such that for any $l\geq l_0$, we have
$$\Hom_{D_{\mathscr{C}}^b(\mathscr{A})}(\Con(s),\Ker d_C^{-l}[l+1]) =0.$$
Take $l=l_0$ in (2). On one hand, applying the functor $\Hom_{D_{\mathscr{C}}^b(\mathscr{A})}(-,\Ker d_C^{-l_0}[l_0+1])$ to (1) we get $h :Z^\bullet\to \Ker d_C^{-l_0}[l_0+1]$
such that $s'=hs$. So we have $s\backslash a=s'\backslash(ha)$. On the other hand, applying $\Hom_{D_{\mathscr{C}}^b(\mathscr{A})}(M,-):=(M,-)$ to (2) we get an exact sequence
$$(M,N)\buildrel {(M,s')} \over \longrightarrow (M,\Ker d_C^{-l_0}[l_0+1])\to (M,(\sigma^{\geq -l_0}C^\bullet)[1]).$$
Since $M\in {^{\bot_\mathscr{C}}\mathscr{C}}$, by using induction on $\omega(\sigma^{\geq -l_0}C^\bullet)$ we have
$(M,(\sigma^{\geq -l_0}C^\bullet)[1])$=0, and hence there exists $f: M\to N$ such that
$ha$=$s'f$. Therefore we have $s\backslash a$=$s'\backslash(h a)$=$s'\backslash(s'f)$=$\id_N\backslash f$, that is, $G$ is surjective.

Next, if $f:M\to N$ satisfies $G(f)=\id_N\backslash f=0$ in $D_{\mathscr{C}{\text-}sg}(\mathscr{A})$, then there exists $s:N\to Z^\bullet$ with $\Con(s)
\in K^b(\mathscr{C})$ such that $sf=0$ in $D^b_{\mathscr{C}}(\mathscr{A})$. Use the same notations as in (1) and (2), by the above argument we have $s'=hs$, so $s'f=0$. Applying
$\Hom_{D_{\mathscr{C}}^b(\mathscr{A})}(M,-)$ to (2) we get that there exists $f^{'}:M\to \sigma^{\geq -l_0}C^\bullet$ such that $f=d_C^0f^{'}$.

Put $\sigma^{<0}(\sigma^{\geq -l_0})C^\bullet:=0\to C^{-l_0}\to C^{-l_0+1}\to\cdots \to C^{-1}\to 0$. We have the following distinguished triangle:
$$(\sigma^{<0}(\sigma^{\geq -l_0})C^\bullet)[-1]\longrightarrow C^0\buildrel {\pi} \over \longrightarrow \sigma^{\geq -l_0}C^\bullet \to\sigma^{<0}(\sigma^{\geq -l_0})C^\bullet$$
in $D_{\mathscr{C}}^b(\mathscr{A})$, where $\pi$ is the canonical map. By applying the functor $\Hom_{D_{\mathscr{C}}^b(\mathscr{A})}(M,-)$
to this triangle, it follows from $M\in {^{\bot_\mathscr{C}}\mathscr{C}}$ that
$\Hom_{D_{\mathscr{C}}^b(\mathscr{A})}(M,\sigma^{<0}(\sigma^{\geq -l_0})C^\bullet)=0$, and hence there exists $g:M\to C^0$ such that $f^{'}=\pi g$. So $f=d_C^0\pi g$ in
$D_{\mathscr{C}}^b(\mathscr{A})$. By Proposition 3.3(3), $\mathscr{A}$ is a full subcategory of $D_{\mathscr{C}}^b(\mathscr{A})$. So $f$ factors through $C^0$ in $\mathscr{A}$,
and hence $\Ker G\subseteq\mathscr{C}(M,N)$. Since $\mathscr{C}(M,N)\subseteq\Ker G$ trivially, $\Ker G=\mathscr{C}(M,N)$,
which means that $\Hom_\mathscr{A}(M,N)/\mathscr{C}(M,N)\cong \Hom_{D_{\mathscr{C}{\text-}sg}(\mathscr{A})}(M,N)$. \hfill{$\square$}

\vspace{0.2cm}

Let $\theta :\mathscr{G}(\mathscr{C})\to D_{\mathscr{C}{\text-}sg}(\mathscr{A})$ be the composition of the following three functors:
 the embedding functors $\mathscr{G}(\mathscr{C})\hookrightarrow\mathscr{A}$,
$\mathscr{A}\hookrightarrow D_{\mathscr{C}}^b(\mathscr{A})$ and the localization functor $D_{\mathscr{C}}^b(\mathscr{A})\to D_{\mathscr{C}{\text-}sg}(\mathscr{A})$, and let
$\underline{\mathscr{G}(\mathscr{C})}$ denote the stable category of $\mathscr{G}(\mathscr{C})$.

\vspace{0.2cm}

{\bf Proposition 4.9.}  {\it $\theta$ induces a fully faithful functor $\theta'$ : $\underline{\mathscr{G}(\mathscr{C})}\to D_{\mathscr{C}{\text-}sg}(\mathscr{A})$.}

\vspace{0.2cm}

{\it Proof.} Since $\mathscr{G}(\mathscr{C})\subseteq {^{\perp_{\mathscr{C}}}\mathscr{C}}$, the assertion follows from Lemma 4.8. \hfill{$\square$}

\vspace{0.2cm}

Recall from [C2] that a {\it $\partial$-functor} is an additive functor $F$ from an exact category $(\mathscr{B},\varepsilon)$ to a triangulated category $\mathcal{C}$ satisfying that
for any short exact sequence $ L\buildrel {i} \over \longrightarrow M\buildrel {p} \over \longrightarrow N$ in $\varepsilon$, there exists a morphism $\omega_{(i,p)}:F(N)\to F(L)[1]$
such that the the triangle
$$F(L)\buildrel {F(i)} \over \longrightarrow F(M)\buildrel {F(p)} \over \longrightarrow F(N)\buildrel {\omega_{(i,p)}} \over \longrightarrow F(L)[1]$$
in $\mathcal{C}$ is distinguished; moreover, the morphism $\omega_{(i,p)}$ are ``functorial" in the sense that any morphism between two short exact sequences in $\varepsilon$:
$$\xymatrix{L\ar[r]^{i}\ar[d]^{f} & M\ar[r]^{p}\ar[d]^{g}& N\ar[d]^h \\
L'\ar[r]^{i'}& M'\ar[r]^{p'}&N',
}$$
the following is a morphism of triangles:
$$\xymatrix{F(L)\ar[r]^{F(i)}\ar[d]^{F(f)} & F(M)\ar[r]^{F(p)}\ar[d]^{F(g)}& F(N)\ar[d]^{F(h)}\ar[r]^{\omega_{(i,p)}}&F(L)[1]\ar[d]^{F(f)[1]} \\
F(L')\ar[r]^{F(i')}& F(M')\ar[r]^{F(p')}&F(N')\ar[r]^{\omega_{(i',p')}}&F(L')[1]
.}$$

\vspace{0.2cm}
By [H1, Chapter I, Theorem 2.6] and Corollary 4.7, $\underline{\mathscr{G}(\mathscr{C})}$ and $D_{\mathscr{C}{\text-}sg}(\mathscr{A})$ are triangulated categories.
Moreover, we have

\vspace{0.2cm}

{\bf Proposition 4.10.} {\it The functor $\theta^{'}$ in Proposition 4.9 is a triangle functor.}

\vspace{0.2cm}

{\it Proof.}  We first claim that $\theta$ is a $\partial$-functor. In fact, let $0\to L\buildrel {f} \over \longrightarrow M\buildrel {g} \over \longrightarrow N\to 0$ be
a $\Hom_\mathscr{A}(\mathscr{C},-)$-exact complex with all terms in $\mathscr{G}(\mathscr{C})$. Then it induces a distinguished triangle in  $D_{\mathscr{C}{\text-}sg}(\mathscr{A})$,
saying $\theta(L)\buildrel {\theta(f)} \over \longrightarrow\theta(M)\buildrel {\theta(g)} \over \longrightarrow\theta(N)\buildrel {\omega_{(f,g)}} \over \longrightarrow\theta(L)[1]$.
It is clear that $\omega_{(f,g)}$ is ``functorial" . This shows that $\theta$ is a $\partial$-functor.

Note that every object in $\mathscr{C}$ is zero in $D_{\mathscr{C}{\text-}sg}(\mathscr{A})$. So $\theta$ vanishes on the projective-injective objects in $\mathscr{G}(\mathscr{C})$.
It follows from [C2, Lemma 2.5] that the induced functor $\theta'$ is a triangle functor.  \hfill{$\square$}

\vspace{0.2cm}

By Propositions 4.9 and 4.10 the natural triangle functor $\underline{\mathscr{G}(\mathscr{C})}\to D_{\mathscr{C}{\text-}sg}(\mathscr{A})$ is fully faithful.
It is of interest to make sense when it is essentially surjective (or dense). We have the following

\vspace{0.2cm}

{\bf Theorem 4.11.} {\it If $\mathscr{C}\mathscr{G}(\mathscr{C})$-dim $\mathscr{A}<\infty$, then the natural functor $\theta :\mathscr{G}(\mathscr{C})\to
D_{\mathscr{C}{\text-}sg}(\mathscr{A})$ is essentially surjective (or dense).}

\vspace{0.2cm}
{\it Proof.} Let $X^\bullet\in D_\mathscr{C}^b(\mathscr{A})$. By Proposition 3.4, there exists $C_0^\bullet=(C_0^i,d_{C_0}^i)\in K^{-,\mathscr{C}b}(\mathscr{C})$
such that $X^\bullet\cong C_0^\bullet$ in $D_\mathscr{C}^b(\mathscr{A})$. So there exists $n_0\in \mathbb{Z}$ such that $H^i(\Hom_{\mathscr{A}}(\mathscr{C},C_0^{\bullet}))=0$ for any
$i \leq n_0$. Let $K^i=\Ker d_{C_0}^i$. Then $C_0^\bullet$ is isomorphic to the complex:
$$0\to K^i\to C_0^i\buildrel {d^i_{C_0}} \over \longrightarrow C_0^{i+1}\buildrel {d^{i+1}_{C_0}} \over \longrightarrow C_0^{i+2}\to\cdots$$
in $D_\mathscr{C}^b(\mathscr{A})$ for any $i \leq n_0$. It induces a distinguished triangle in $D_\mathscr{C}^b(\mathscr{A})$, hence a distinguished triangle in
$D_{\mathscr{C}{\text-}sg}(\mathscr{A})$ of the following form:
$$K^i[-i]\to \sigma^{\geq i}C_0^\bullet\to C_0^\bullet\to K^i[-i+1].$$  Since
$\sigma^{\geq i}C_0^\bullet\in K^b(\mathscr{C})$, $C_0^\bullet\cong K^i[-i+1]$ in $D_{\mathscr{C}{\text-}sg}(\mathscr{A})$. Take $l_0=i$ and $Y=K^i$.
Then $C_0^\bullet\cong Y[-l_0+1]$ in $D_{\mathscr{C}{\text-}sg}(\mathscr{A})$. By assumption we may assume that $\mathscr{C}\mathscr{G}(\mathscr{C})$-$\dim Y=m_0<\infty$.
Let $C_1^\bullet\to Y$ be the left $\mathscr{C}$-resolution of $Y$.
We claim that for any $n\leq-m_0+1$, $\Ker d_{C_1}^{n}\in\mathscr{G}(\mathscr{C})$, where $d_{C_1}^{n}$ is the $n$th differential of $C_1^\bullet$.

We have a $\mathscr{C}$-acyclic complex:
$$0\to G^{-m_0}\to G^{-m_0+1}\to \cdots\to G^{-1}\to G^{0}\to Y\to 0$$
with $G^{j}\in \mathscr{G}(\mathscr{C})$ for any $-m_0\leq j\leq0$. Let $G^\bullet$ be the complex $0\to G^{-m_0}\to G^{-m_0+1}\to \cdots\to G^{-1}\to G^{0}\to 0$.
By Lemma 2.3, there exists a $\mathscr{C}$-quasi-isomorphism $C_1^\bullet\to G^\bullet$ lying over $\id_Y$, and hence its mapping cone is $\mathscr{C}$-acyclic.
So for any $n\leq-m_0+1$, we get the following $\mathscr{C}$-acyclic complex:
$$0\to \Ker d_{C_1}^{n}\to C_1^n\to \cdots\to C_1^{-m_0}\to C_1^{-m_0+1}\oplus G^{-m_0}\to \cdots\to C_1^0\oplus G^{-1}\to G^0\to 0.$$
Note that this complex is acyclic because $\mathscr{C}$ is admissible. Put $K=\Ker(C_1^0\oplus G^{-1}\to G^0)$, we get a $\Hom_{\mathscr{A}}(\mathscr{C},-)$-exact exact sequence
$0\to K\to C_1^0\oplus G^{-1}\to G^0\to 0$. By Lemma 3.7(3), we get an exact sequence:
$$0\to \Hom_\mathscr{A}(G^0,C)\to \Hom_\mathscr{A}(C_1^0\oplus G^{-1},C)\to \Hom_\mathscr{A}(K,C)\to \Ext_\mathscr{C}^1(G^0,C)$$
for any $C\in \mathscr{C}$. Since $G^0\in\mathscr{G}(\mathscr{C})$, $\Ext_\mathscr{C}^1(G^0,C)=0$ and so the exact sequence
$0\to K\to C_1^0\oplus G^{-1}\to G^0\to 0$ is $\Hom_{\mathscr{A}}(-,\mathscr{C})$-exact. Because both $C_1^0\oplus G^{-1}$ and $G^0$ are in $\mathscr{G}(\mathscr{C})$,
$K\in \mathscr{G}(\mathscr{C})$ by [Hu, Proposition 4.7]. Iterating this process, we get that $\Ker d_{C_1}^{n}\in \mathscr{G}(\mathscr{C})$ for any $n\leq-m_0+1$.
The claim is proved.

Choose a left $\mathscr{C}$-resolution $C_1^\bullet$ of $Y$ and put $X=\Ker d_{C_1}^{-m_0+1}$. By the above claim we have a $\mathscr{C}$-acyclic complex:
$$0\to X\to C_1^{-m_0+1}\to C_1^{-m_0+2}\to \cdots\to C_1^{0}\to Y\to 0$$
with $X\in \mathscr{G}(\mathscr{C})$. Then $Y\cong X[m_0]$ in $D_{\mathscr{C}{\text-}sg}(\mathscr{A})$ and $X^\bullet\cong C_0^\bullet\cong Y[-l_0+1]\cong X[m_0-l_0+1]$
in $D_{\mathscr{C}{\text-}sg}(\mathscr{A})$. We may assume that $X^\bullet\cong C_0^\bullet\cong X[r_0]$ in $D_{\mathscr{C}{\text-}sg}(\mathscr{A})$ for $r_0>0$.
Because $X\in \mathscr{G}(\mathscr{C})$, we get a $\Hom_{\mathscr{A}}(\mathscr{C},-)$-exact exact sequence $0\to X\to C^0\to C^1\to\cdots\to C^{r_0-1}\to X'\to 0$
with $X'\in\mathscr{G}(\mathscr{C})$ and $C^i\in \mathscr{C}$ for any $0\leq i\leq r_0-1$. It follows that $X\cong X'[-r_0]$ and
$X^\bullet \cong C_0^\bullet\cong X[r_0]\cong X'$ in $D_{\mathscr{C}{\text-}sg}(\mathscr{A})$. This completes the proof. \hfill{$\square$}

\vspace{0.2cm}

The following is the main result of this paper.

\vspace{0.2cm}

{\bf Theorem 4.12.} {\it If $\mathscr{C}\mathscr{G}(\mathscr{C})$-$\dim \mathscr{A}<\infty$, then the natural functor
$\theta: \mathscr{G}(\mathscr{C})\to D_{\mathscr{C}{\text-}sg}(\mathscr{A})$ induces a triangle-equivalence
$\theta^{'}: \underline{\mathscr{G}(\mathscr{C})} \to D_{\mathscr{C}{\text-}sg}(\mathscr{A})$.}

\vspace{0.2cm}
{\it Proof.} It follows directly from Propositions 4.9, 4.10 and Theorem 4.11. \hfill{$\square$}

\vspace{0.2cm}

The following result is the dual version of Happel's result, see [H2, Theorem 4.6].

\vspace{0.2cm}

{\bf Corollary 4.13.} {\it If $A$ is Gorenstein, then the canonical functor $\mathscr{G}(A)\to D_{sg}(A) $ induces a triangle-equivalence
$\underline{\mathscr{G}(A)}\to D_{sg}(A)$.}

\vspace{0.2cm}

{\it Proof.} Let $A$ be Gorenstein and $\mathscr{C}=A{\text -}\proj$. Then $\mathscr{C}\mathscr{G}(\mathscr{C})$-$\dim \mathscr{A}<\infty$
by [Hos, Theorem]. Now the assertion is an immediate consequence of Theorem 4.12. \hfill{$\square$}

\vspace{0.5cm}

{\bf Acknowledgements.} This research was partially supported by NSFC (Grant No. 11171142)
and a Project Funded by the
Priority Academic Program Development of Jiangsu Higher Education Institutions. The authors thank the referee for the helpful suggestions.

\end{document}